\documentclass[11pt,oneside,notitlepage]{amsart}
\usepackage[german]{babel}
\usepackage[latin1]{inputenc}
\usepackage{mathpazo}
\usepackage{footmisc}

\vfuzz \hfuzz
\let\Person\textsc
\let\Emphasis\textbf
\let\Author\textsc
\def\Editor{}
\let\Title\textit
\let\Place\textit
\let\Foreign\textit

\begin{document}

\begin{center}
\textsc{\Large\bf Ueber die Geometrie der alten \AE{}gypter}\footnote{Presented
to the Royal Academy at Vienna. Originally published by the Royal Press
»in Commission bei Karl Gerold's Sohn«, 1884.
Transcribed by Ralf Stephan,
eMail: \texttt{mailto:ralf@ark.in-berlin.de}.
Other formats at \texttt{http://www.gutenberg.org/etext/24817}}\\[5pt]
{\Small Vortrag, gehalten in der feierlichen Sitzung der Kaiserlichen Akademie der Wissenschaften am 29. Mai 1884.}\\[10pt]
Dr. Emil Weyr\\[15pt]
\end{center}

\bigskip

MSC-Class: 51-03 01A16

\bigskip

Möge mir gestattet sein, bei dem heutigen feierlichen Anlasse
ein Bild zu entrollen, welches in grossen Strichen die
allgemeinen Umrisse des Zustandes der geometrischen
Wissenschaften bei den alten Aegyptern zur Darstellung
bringen soll; und möge dasselbe Wohlwollen, das, gepaart
mit einer althergebrachten Sitte, mich heute auf diesen eben so
ehrenvollen als schwierigen Platz gestellt, auch bei der
Beurtheilung der folgenden bescheidenen, weil schwachen
Kräften entspringenden Leistung obwalten!

So wie der Anfang aller menschlichen Kenntnisse, so ist
auch der Ursprung der Geometrie in grauestes Alterthum zu
versetzen, er ist zu suchen in jenen der Zeit nach unangebbaren
Perioden der menschlichen Entwicklung, in welchen
das erste Erwachen des Selbstbewusstseins zu finden wäre.
Sind doch manche geometrische Anschauungen auch dem
Thiere eigen; so jene der geraden Verbindungslinie zweier
Punkte als der kürzesten Entfernung; jene des Mehr und
Weniger bei Quantitäten der Entfernungen, Höhen, Neigungen,
und so werden auch manche abstractere Raumanschauungen
dem Menschen in seinen ersten Entwicklungsperioden eigen
geworden sein, Anschauungen, welche durch die Möglichkeit
und auf Grund der sprachlichen Bezeichnung jene Stabilität
erhielten, die sie befähigte, als erste Fundamente der geometrischen
Kenntnisse zunächst, und der Geometrie als Wissenschaft
später aufzutreten.


Geometrisches Denken entstand zu den verschiedensten
Zeiten, an den verschiedensten Orten. Denn überall, wo
der menschliche Geist sich zu entwickeln begann, und
das menschliche Denken jene Höhe erreichte, auf welcher
Abstractionen entstehen, bildeten sich die grundlegenden
Raumbegriffe; der des Punktes, der geraden und krummen
Linien, der ebenen und krummen Flächen. Denn überall
in der Natur boten sich dem erwachenden Menschen
Repräsentanten dieser Begriffe in grösserer oder geringerer
Genauigkeit dar. Während der Anblick der auf- und untergehenden
Sonne, sowie des vollen Mondes in südlichen
Gegenden fast täglich das Bild der »vollkommensten«,
der »schönsten« Linie, der Kreislinie vorführte, stellten sich
die zahllosen Sterne des Abends dem Auge als glänzende
Punkte dar, welche in ihren mannigfaltigen gegenseitigen
Lagenverhältnissen die Phantasie des Menschen bei der, von
ihm beliebten Eintheilung des Himmels in Sternbilder zur
Herstellung so mancher geraden und krummen Linien verleiten
mochten. Und selbst in seiner nächsten Umgebung
fand der beobachtende Mensch geometrische Anklänge; das
Gewebe der Spinne mit seinen kreisrunden und radialen
Fäden, die sechseckige Bienenzelle, die beim Fallen eines
Körpers in ruhendes Wasser entstehenden concentrischen
Wellenringe, und wie vieles Andere musste, wenn auch
nach und nach, so doch mit zwingender Nothwendigkeit
den Menschen zur Beobachtung gesetzmässiger geometrischer
Formen führen.

Als Mutterland der Mathematik im Allgemeinen, und
der Geometrie im Besonderen wird Aegypten angeführt; doch
ist die Zeit längst vorbei, wo man sich Aegypten als einzigen
Ursprungsort dieser Wissenschaften dachte, vielmehr muss
als feststehend angenommen werden, dass jedes Volk in
seinem Entwicklungsgange geometrische Anschauungen sich
anzueignen schon durch praktische Bedürfnisse gezwungen
war. Die Höhe, zu welcher sich die einzelnen Völker in ihren
mathematischen Speculationen emporzuschwingen vermochten,
hing von der Richtung des Bildungsganges, von
dem Maasse des Bedürfnisses und nicht in letzter Reihe von
dem Einflüsse religiöser Verhältnisse ab.

Und so mag sich zunächst jene Naturgeometrie entwickelt
haben, welche allen Völkern zugesprochen werden
muss, und auf deren Vorhandensein, weil auf die Anwendungen
ihrer freilich einfachsten Principien, Ueberreste von
Bauten überall dort hinweisen, wo wir in der Lage sind,
solche beobachten zu können. Die Pellasger, die vorhellenischen
Ureinwohner Griechenlands, mussten lange vor Entstehung
der Philosophie geometrische Kenntnisse in dem
Maasse besessen haben, wie sie zur Aufführung von Wasserbauten,
Dämmen, Canälen und Burgen, von denen man jetzt
noch Spuren findet, nothwendig waren.

Verfolgt man die Entwicklung der Geometrie zu ihren
Quellen aufwärts, so dürfen wir nicht überrascht sein, dass
man bei dem uns bekannten ältesten Culturvolke, bei den
Aegyptern, am weitesten vorzudringen vermag, und zwar an
der Hand der indirecten wie der directen Nachrichten, welche
uns über diesen Gegenstand zugekommen sind. Leider jedoch
sind die Ersteren ihrem Inhalte und die Letzteren ihrer Zahl
nach nur spärliche zu nennen.

Zahlreich sind wohl die Stellen in griechischen Philosophen
und Geschichtschreibern, welche Bezug haben auf
aegyptische Geometrie, es lässt sich jedoch nicht verkennen,
dass oft die Späteren auf Frühere sich stützen, und wir es
möglicherweise mit einer einzigen, durch Jahrhunderte fortgeführten
Nachricht zu thun haben.


  Durch \Person{Herodot}, welcher um die Mitte des fünften
  vorchristlichen Jahrhunderts (460) Aegypten bereiste, erfahren
wir\footnote{\Author{Herodot}, \Title{Reisebericht},
    II, 109.},
dass die Geometrie von Aegypten nach Griechenland
verpflanzt worden sei. Etwas später (393 v.~Chr.)
berichtet \Person{Isokrates} die Thatsache\footnote{
  \Author{Isokrates}, \Title{Busiris},
    c. 9.}, dass die Aegypter »die
Aelteren (unter ihren Priestern) über die wichtigsten Angelegenheiten
setzten, dagegen die Jüngeren beredeten, mit
Hintansetzung des Vergnügens, sich mit Astronomie, Rechenkunst
und Geometrie zu beschäftigen«.

In \Person{Platon}'s \Title{Phädrus} sagt
\Person{Sokrates}: »Ich habe
vernommen, zu Naukratis in Aegypten sei einer der dortigen
alten Götter gewesen, dem auch der Vogel geheiligt ist,
den sie Isis nennen, während der Gott selbst den Namen
Teuth führt; dieser habe zuerst Zahlenlehre und Rechenkunst
erfunden und Geometrie und Astronomie«\footnote{
\Title{\Emphasis{Platonis} Phaedrus},
  ed. \Editor{Ast.} 
  I. p. 246.},
und einen
directen Hinweis finden wir bei \Person{Aristoteles}, welcher in
seiner \Title{Metaphysik} sagt:\footnote{
  \Author{Aristoteles}, \Title{Metaph. I},
    1.}
»Daher entstanden auch in Aegypten
die mathematischen Wissenschaften, denn hier war den
Priestern die dazu nöthige Müsse vergönnt.«

Uebrigens schrieben sich die Aegypter neben der Erfindung
der Buchstabenschrift auch jene der meisten Wissenschaften
und Künste zu, worüber \Person{Diodor}\footnote{
  \Author{Diodor}, I, 69.},
welcher etwa
70 Jahre v.~Chr.~G. Aegypten bereiste, bemerkt: »Die
Aegypter behaupten, von ihnen sei die Erfindung der Buchstabenschrift
und die Beobachtung der Gestirne ausgegangen,
ebenso seien von ihnen die Theoreme der Geometrie und
die meisten Wissenschaften und Künste erfunden worden.«

Neben diesen ganz allgemein gehaltenen Angaben sind
hauptsächlich diejenigen Berichte zu erwähnen, welche sich
auf die Art der wissenschaftlichen Leistungen der Aegypter
beziehen.


  Da sagt zunächst \Person{Herodot}\footnote{
    Herodot l. c.}
in Hinsicht auf die unter
dem Könige \Person{Sesostris} durchgeführte Ländereintheilung:
»Auch sagten sie, dass dieser König das Land unter alle
Aegypter so vertheilt habe, dass er jedem ein gleich grosses
Viereck gegeben, und von diesem seine Einkünfte bezogen habe,
indem er eine jährlich zu entrichtende Steuer auflegte.
Wem aber der Fluss (Nil) von seinem Theile etwas
wegriss, der musste zu ihm kommen und das Geschehene
anzeigen; er schickte dann die Aufseher, die auszumessen
hatten, um wie viel das Landstück kleiner geworden war,
damit der Inhaber von dem übrigen nach Verhältniss der
aufgelegten Abgaben steure. Hieraus erscheint mir die Geometrie
entstanden zu sein, die von da nach Hellas kam.«

Die, \Person{Herodot}, dem Vater der Geschichtsschreibung
folgenden Berichterstatter hielten sich nun, vielleicht erklärlicherweise,
vorzüglich an den einen, die Nilüberschwemmungen
betreffenden Theil obiger Nachricht, und wurde,
gewiss Unberechtigtermassen der Nil als der unmittelbare
Anstoss für alle geometrischen Arbeiten der Aegypter hingestellt.
Und doch scheint es uns viel näherliegend, die einerseits
behufs der Steuerbemessung und Controle, anderseits
wegen der aus den Veränderungen im Besitzstande sich
nothwendig ergebenden Flächenfestsetzungen als den Hauptbeweggrund
jener Vermessungen zu erkennen, wobei die
gesammelten Erfahrungen gewiss auch bei der Beurtheilung
der unzweifelhaft nach den periodisch eintretenden Nilüberschwemmungen
vorgekommenen Terrainveränderungen mit
Vortheil benutzt worden sein mögen.

Unverkennbar ist der Zug nach Aufbauschung und Ausschmückung
des, jene Nilüberschwemmungen betreffenden
Theiles des \Person{Herodot}'schen Berichtes, wenn man die Aufzeichnungen
späterer Gewährsmänner näher betrachtet.


  Zunächst finden wir bei \Person{Heron} dem Aelteren die folgende
  diesbezügliche Stelle\footnote{
    \Title{\Emphasis{Heronis Alexandr.} geom. et stereom. reliquiae},
    ed. \Editor{Hultsch.}
    p. 138.}:
  »Die früheste Geometrie beschäftigte
sich, wie uns die alte Ueberlieferung lehrt, mit der Messung und
Vertheilung der Ländereien, woher sie Feldmessung genannt
wurde. Der Gedanke einer Messung nämlich ward den Aegyptern
an die Hand gegeben durch die Ueberschwemmungen des
Nil. Denn viele Grundstücke, die vor der Flussschwelle offen
dalagen, verschwanden beim Steigen des Flusses und kamen
erst nach dem Sinken desselben zum Vorschein, und es war
nicht immer möglich, über die Identität derselben zu entscheiden.
Dadurch kamen die Aegypter auf den Gedanken
einer solchen Messung des vom Nil blossgelegten Landes.«

Weiter finden wir bei \Person{Diodor}\footnote{
  \Author{Diodor},  I, 81.}
einen Ausspruch,
durch welchen wir übrigens auch über andere wissenschaftliche
Leistungen der Aegypter belehrt werden; \Person{Diodor}
sagt: »Die Priester lehren ihre Söhne zweierlei Schrift, die
sogenannte heilige, und die, welche man gewöhnlich lernt.
Mit Geometrie und Arithmetik beschäftigen sie sich eifrig.
Denn indem der Fluss jährlich das Land vielfach verändert,
veranlasst er viele und mannigfache Streitigkeiten über die
Grenzen zwischen den Nachbarn; diese können nun nicht
leicht ausgeglichen werden, wenn nicht ein Geometer den
wahren Sachverhalt durch directe Messung ermittelt. Die
Arithmetik dient ihnen in Haushaltungsangelegenheiten und
bei den Lehrsätzen der Geometrie; auch ist sie denen von
nicht geringem Vortheile, die sich mit Sternkunde beschäftigen.
Denn wenn bei irgend einem Volke die Stellungen und Bewegungen
der Gestirne sorgfältig beobachtet worden sind,
so ist es bei den Aegyptern geschehen; sie verwahren Aufzeichnungen
der einzelnen Beobachtungen seit einer unglaublich
langen Beihe von Jahren, da bei ihnen seit alten
Zeiten her die grösste Sorgfalt hierauf verwendet worden
ist. Die Bewegungen und Umlaufszeiten sowie die Stillstände
der Planeten, auch den Einfluss eines jeden auf die Entstehung
lebender Wesen und alle ihre guten und schädlichen
Einwirkungen haben sie sehr sorgfältig beobachtet.«

Am innigsten verknüpft erscheint die Geometrie der
Aegypter mit den Ueberschwemmungen des Nil bei \Person{Strabon}\footnote{
  \Author{Strabon}, ed. \Editor{Meinike},
    lib. XVII, C. 787, p. 1098.};
welcher bemerkt, »dass es einer sorgfältigen und bis auf das
Genaueste gehenden Eintheilung bedurfte, wegen der beständigen
Verwüstung der Grenzen, die der Nil bei seinen
Ueberschwemmungen veranlasst, indem er Land wegnimmt
und zusetzt, und die Gestalt verändert, und die anderen Zeichen
unkenntlich macht, wodurch das fremde und eigene Besitzthum
unterschieden wird. Man müsse daher immer und
immer wieder messen. Hieraus soll die Geometrie entstanden
sein.«

Den gesellschaftlichen Einrichtungen der Aegypter
entsprechend, muss als feststehend angenommen werden,
dass sich eine Kaste, nach eben Gehörtem die der Priester,
mit dem wissenschaftlichen Theile der Geometrie beschäftigte,
während eine andere, die der Feldmesser, die von den
Ersteren aufgestellten und sorgsam gehüteten geometrischen
Principien praktisch zur Anwendung brachte. Dabei wurden,
wie wir später sehen werden, die Geheimnisse der Priester,
insoweit sie geometrische Wahrheiten und Berechnungsregeln
betrafen, möglicherweise nur insoweit enthüllt, dass
bei deren Verwendung nur annäherungsweise richtige Resultate
zum Vorschein kamen.

Wohl sind einige Schriftsteller so weit gegangen, dass
sie, die unläugbaren Uebertreibungen des Zusammenhanges
zwischen den Nilüberschwemmungen und der ägyptischen
Geometrie im Auge behaltend, die Existenz der letzteren
einfach negirten, und alle die citirten Aussprüche in das
Gebiet der Fabel verwiesen.

Was macht man jedoch dann mit den wohlbeglaubigten
Nachrichten über die Reisen, welche hervorragende griechische
Philosophen nach Aegypten unternahmen, oft jahrelang
dort verweilend, um sich in die Geheimnisse aegyptischer
Priester einweihen und mit deren geometrischem Wissen
vertraut machen zu lassen?

\Person{Eudemus} von \Place{Rhodos}\footnote{
  \Title{\Emphasis{Eudemi Rhodii} Peripatetici fragmenta quae
      supersunt}.
    ed. \Editor{L. Spengel.}
    Berlin 1870.},
  einer der ältesten Peripatetiker,
schrieb eine Geschichte der Mathematik, aus
welcher uns durch \Person{Proklos Diadochus}\footnote{
  \Title{\Emphasis{Procl.} comment.}
    ed. \Editor{Rasil.}
    p. 19;
  \Author{Barocius}
    p. 37.}, einen Philosophen
des fünften nachchristlichen Jahrhunderts, ein
Bruchstück erhalten ist, welches sozusagen das einzige
Mittel bildet, das uns einen Einblick in die geometrischen
Errungenschaften der Griechen in den ersten dritthalb
Jahrhunderten nach \Person{Thales} gewährt. Hierin heisst es unter
Anderem: »\Person{Thales}, der nach Aegypten ging, brachte zuerst
die Geometrie nach Hellas hinüber und Vieles entdeckte er
selbst, von Vielem aber überlieferte er die Anfänge seinen
Nachfolgern; das Eine machte er allgemeiner, das Andere
mehr sinnlich fassbar.« Hundert Jahre nach dem Tode des
\Person{Pythagoras} berichtet der Redner \Person{Isokrates}\footnote{
  \Author{Isokrates}, \Title{Busiris},
    cap. 11.}: »Man
könnte, wenn man nicht eilen wollte, viel Bewunderungswürdiges
von der Heiligkeit aegyptischer Priester anführen,
welche ich weder allein noch zuerst erkannt habe, sondern
viele der jetzt Lebenden und der Früheren, unter denen
auch \Person{Pythagoras} der Samier ist, der nach Aegypten
kam und ihr Schüler wurde und die fremde Philosophie
zuerst zu den Griechen verpflanzte.«

Während der Aufenthalt des \Person{Pythagoras} in Aegypten
unter Anderen auch noch von \Person{Strabon}\footnote{
  \Author{Strabon},
    XIV, 1. 16.}
und \Person{Antiphon}\footnote{
  \Author{Porphyrius}, \Title{De vita Pythagorae}
    cap. 7;
  \Author{Diogenes Laertius},
    VIII, 3.}
bestätiget wird, nennt uns \Person{Diodor}\footnote{
  \Author{Diodor}, I, c. 96.}
eine ganze Reihe von
Namen, indem er sagt; »Die aegyptischen Priester nennen
unter den Fremden, welche nach den Verzeichnissen in den
heiligen Büchern vormals zu ihnen gekommen seien, den
\Person{Orpheus}, \Person{Musaios},
\Person{Melampus} und \Person{Daidalos}, nach
diesen den Dichter \Person{Homer}, den Spartaner \Person{Lykurgos}, ingleichen
den Athener \Person{Solon} und den Philosophen \Person{Platon}.
Gekommen sei zu ihnen auch der Samier \Person{Pythagoras} und
der Mathematiker \Person{Eudoxos}, ingleichen \Person{Demokritos} von
\Place{Abdera} und \Person{Oinopides} von \Place{Chios}. Von allen diesen
weisen sie noch Spuren auf, von den Einen Bildnisse von
den Anderen Orte und Gebäude, die nach ihnen benannt
sind. Aus der Vergleichung dessen, was jeder von ihnen in
seinem Fache geleistet hat, führen sie den Beweis, dass sie
Dasjenige um desswillen sie von den Hellenen bewundert
werden, aus Aegypten entlehnt haben.« Aus diesen Stellen
geht mit Sicherheit hervor, dass viele Griechen nach
Aegypten zogen, um bei den dortigen Priestern Philosophie
und Mathematik kennen zu lernen, da wohl in den Berichten
nur die hervorragenden Männer angeführt wurden.

Der Milesier \Person{Thales}, welcher erst in vorgerücktem
Alter, und nachdem er als Handelsmann früher gewiss schon
mehrmals Aegypten besucht gehabt, sich daselbst behufs
seiner Studien zu längerem Aufenthalt niederlies, ist merkwürdiger
Weise in dem Berichte des Diodor nicht angeführt,
und könnte man wohl aus diesem Umstande umsomehr
einen gewissen Grad von Unglaublichkeit ableiten, als darin
mythische Namen wie \Person{Orpheus}, \Person{Daidalos} und \Person{Homer}
angeführt erscheinen. Diese letzteren konnten jedoch sehr
wohl dem im Ganzen und Grossen sonst richtigen Verzeichnisse
vom Berichterstatter eigenwillig beigefügt worden sein,
um dadurch das hohe Alter aegyptischer Wissenschaft in ein
vorteilhaftes Licht zu setzen.


Abgesehen jedoch von aller Wahrscheinlichkeit oder
Unwahrscheinlichkeit für die Exactheit obiger Aussprüche
in Bezug auf einzelne Namen, dürfte jedenfalls das als
unumstössliche Wahrheit gelten, dass die ägyptischen Priester
von den Griechen als in den Wissenschaften, insbesondere
in der Geometrie sehr bewandert gehalten wurden, und
zwar in einem solchen Maasse, dass eine Reihe hervorragender
griechischer Philosophen es nicht verschmähte, die, für
damalige Verhältnisse nicht unbedeutende Reise nach
Aegypten zu unternehmen, ja oft jahrelang in diesem Lande
mit unbekannter Sprache und Schrift zu verweilen, um sich
die Kenntnisse der Aegypter anzueignen.

Stellt man nun zunächst die Frage nach Quantität und
Qualität des geometrischen Wissens, welches die Griechen
von ihren Studienreisen mit nach Hause brachten, so scheint
dies, selbst vom Standpunkte der unmittelbar nachpythagoräischen
Geometrie, äusserst Weniges gewesen zu sein.

\Person{Thales} von Milet, einer der sieben griechischen
Weltweisen, der Begründer der ionischen Schule, \Person{Thales},
welcher für das Jahr 585 v.~Chr.~G. eine, auch eingetroffene
Sonnenfinsterniss vorherzusagen wusste, soll, den uns von
\Person{Proklos} zugekommenen Berichten zufolge, in Aegypten
nicht viel mehr erfahren haben, als die Sätze über die Gleichheit
der Winkel an der Basis eines gleichschenkligen Dreieckes,
die Gleichheit der Scheitelwinkel am Durchschnitt
zweier Geraden; er wusste ferner, wie ein Dreieck durch eine
Seite und die beiden anliegenden Winkel bestimmt erscheint,
diese Erörterung zur Messung der Entfernungen von Schiffen
auf dem Meere benützend, es war ihm bekannt, dass ein
Kreis durch einen Durchmesser halbirt wird,\footnote{
  \Author{Proklos}, ed. \Editor{Friedlein},
    250, 299, 352, 157.}
und soll er die
Höhe der Pyramiden aus der Länge des Schattens gemessen
haben, höchst wahrscheinlich in dem Momente, wo die
Schattenlänge eines senkrechten Stabes der Stablänge gleich
ist,\footnote{
  \Author{Diogenes Laertius},
    I, 27.
  \Author{Plinius}, \Title{Hist. nat.}
    XXXVI, 12, 17.}
möglicherweise jedoch, wie \Person{Plutarch}\footnote{
  \Author{Plutarch}, ed. \Editor{Didot.}
    Vol. 2, III, p. 174.}
berichtet,
auch zu einer beliebigen Tageszeit. Auch wird ihm von
\Person{Pamphile}\footnote{
  \Author{Diogenes Laertius}
    I, 24--25.}
die Kenntniss des Satzes zugeschrieben, dass
der Peripheriewinkel im Halbkreise ein rechter sei. Gewiss
hat Thales wenigstens jene geometrischen Fundamente in
Aegypten kennen gelernt, welche es ihm ermöglichten, die
genannten Sätze als wahr zu erkennen, wenn auch bei ihm,
selbst bei diesen einfachen Dingen an einen strengen Beweis
nicht gedacht werden kann.

Es wäre jedoch voreilig, aus der Geringfügigkeit der
Thaletischen geometrischen Kenntnisse mit \Person{Montucla}
\footnote{\Author{Montucla}, \Title{Hist. d. math.}
    2. édit. t. I, p. 49.}
zu schliessen, dass auch die Aegypter nicht viel mehr gewusst
hätten. Man kann wohl annehmen, dass die aegyptischen
Priester bei ihrer den Fremden gegenüber beobachteten
Zurückhaltung nur einen Theil ihres Wissens offenbarten;
wer könnte jedoch bemessen, in welchem Verhältnisse dieser
Theil zu ihrem Gesammtwissen stand? Der Ansicht
\Person{Montucla}'s kann man entgegensetzen, dass die Aegypter
den Fremden nur einen kleinen Bruchtheil ihres sorgsam
im Verborgenen gehüteten Wissens preisgegeben haben
mochten, wobei ferner nicht unberücksichtigt bleiben darf,
dass den nach Aegypten gekommenen Griechen auch die
Unkenntniss der Sprache und der Schrift weitere, nicht zu
unterschätzende Schwierigkeiten bereitete, in dem Maasse als
vielleicht Manches, was ihnen die aegyptischen Priester von
aegyptischem Wissen zur Verfügung stellten, unverstanden
bleiben konnte.

Was nun das Wesen aegyptischer Geometrie betrifft, so
finden wir in den Berichten der Alten fast gar keine Anhaltspunkte,
um uns hierüber Klarheit verschaffen zu können, und
war man bis vor Kurzem darauf hingewiesen, aus den
Anfängen griechischer Mathematik auf den Stand der aegyptischen
zurückzuschliessen, was, wie aus dem Vorhergesagten
folgen dürfte, mit nicht geringen Schwierigkeiten verbunden
erscheint.

Die Ansicht, dass die Geometrie der Aegypter eigentlich
nur constructiver Natur war, ähnlich dem was wir als Reisskunst
zu bezeichnen pflegen,\footnote{
  \Author{Bretschneider},
    \Title{Die Geometrie und die Geometer vor Euklides},
    p. 11.
Dem Werke Bretschneiders, sowie jenem \Author{Cantor}'s: 
  \Title{Vorlesungen über Geschichte
    der Mathematik}, sind die grundlegenden Gedanken entnommen.}
dürfte sich nicht als stichhältig
erweisen; es möge jedoch gleich jetzt darauf hingedeutet
werden, dass die Aegypter im Construiren geometrischer
Formen nicht unbewandert sein konnten.

So sagt in etwas prahlerischer Weise \Person{Demokritos}
von \Place{Abdera}\footnote{
  \Author{Clemens Alexandrinus}, \Title{Stromata},
    ed. \Editor{Potter},
    I, 357.}
um 420 v.~Chr.~G.: »Im Construiren von
Linien nach Maassgabe der aus den Voraussetzungen zu
ziehenden Schlüsse hat mich keiner je übertroffen, selbst
nicht die sogenannten Harpedonapten der Aegypter«; und
\Person{Theon} von \Place{Smyrna}\footnote{
  \Author{Theon Smyrnaios}, \Title{lib. de astron.}
  ed. \Editor{Martin},
  p. 272.}
erzählt, dass
»Babylonier, Chaldäer
und Aegypter eifrig nach allerhand Grundgesetzen und
Hypothesen suchten, durch welche den Erscheinungen
genügt werden könnte; zu erreichen suchten sie dies dadurch,
dass sie das früher Gefundene in Ueberlegung zogen, und
über die zukünftigen Erscheinungen Vermuthungen aufstellten,
wobei die Einen sich arithmetischer Methoden bedienten,
wie die Chaldäer, die Anderen construirender wie
die Aegypter«.

Aus diesen und ähnlichen Berichten, sowie aus dem
Umstande, dass die Anfänge der griechischen Geometrie
selbst hauptsächlich constructiver Natur waren, muss man
zu dem Schlusse kommen, dass die alten Aegypter seit unvordenklichen
Zeiten die Reisskunst pflegten, und in der langen
Reihe der Jahrhunderte sicherlich eine ziemlich bedeutende
Masse sowohl einfacher als complicirterer Constructionen
erfanden und in ein gewisses System brachten, von Ersteren
zu Letzteren aufsteigend. Diese Constructionen dürften ihrem
grösseren Theile nach, und zwar jenem Theile nach, welcher,
wenn auch ohne Begründung Gemeingut der die Künste und
Gewerbe betreibenden Kasten wurde, nur solche gewesen
sein, die dem praktischen Bedürfnisse dienen konnten, also
zumeist Ornamentenconstructionen. Wir bemerken hier unter
Anderem das Vorkommen regelmässiger geometrischer
Figuren auf uralten Wandgemälden, wie sie sich z.~B. als
färbige Zeichnungen aus den Zeiten der fünften Dynastie,
also unmittelbar nach den Erbauern der Pyramiden, das ist
3400 Jahre v.~Chr.~G. etwa vorfinden.\footnote{
  \Author{Prisse d'Avennes},
    \Title{Hist. de l'art Egypt. d'après les monuments.}}

Man sieht unter der grossen Menge der in dieser Zeit
vorkommenden Figuren eine, aus verschobenen, ineinander
gezeichneten, theilweise durch zu einer Diagonale Parallele
zerlegten Quadraten zusammengesetzte Figur, ferner aus der
Zeit von der zwölften bis zur sechsundzwanzigsten Dynastie,
eine Figur, bestehend aus einem Quadrate, und zwei, längs
der Diagonale centrisch hineingelegten lemniscatischen
Curven, sowie eine Zusammenstellung von um fünfundvierzig
Grade gegeneinander verdrehten, sich durchsetzenden Quadraten.
Kreise erscheinen durch ihre Durchmesser in gleiche
Kreisausschnitte getheilt; so zunächst durch zwei oder vier
Durchmesser in vier beziehungsweise acht, und in späteren
Zeiten auch durch sechs Durchmesser in zwölf gleiche Ausschnitte;
die in den Zeichnungen vorkommenden Wagenräder
besitzen zumeist sechs, seltener vier Speichen, so dass
auch die Theilung des Kreises durch drei Diameter in sechs
gleiche Kreisausschnitte vertreten erscheint.

In einer unvollendet gebliebenen Kammer des Grabes
\Person{Seti~I.}, des Vater
\Person{Ramses~II.} aus der neunzehnten Dynastie
(das sogenannte Grab \Person{Belzoni})\footnote{
  \Author{Wilkinson},
    \Title{Manners and customs of the ancient Egyptians},
    III, p. 313.}
finden wir die Wände
behufs Anbringung von Reliefarbeiten mit einem Netze gleich
grosser Quadrate bedeckt, und es kann keinem Zweifel unterliegen,
dass wir es hier mit der Anwendung eines Verkleinerungs-
beziehungsweise Vergrösserungsmaassstabes zu thun
haben.

Wenn nun auch die einfachen Figuren des Dreieckes,
Quadrates und des Kreises höchst wahrscheinlich ohne
besondere Ueberlegung, einfach dem inneren geometrischen
Formendrange entsprungen sein dürften, so ist doch gewiss,
dass ihre verschiedenartige Zusammensetzung zu Mustern
das Product, wenn auch primitiven geometrischen Denkens
war, welches dann schon eine ziemliche Selbstständigkeit
erreicht haben musste, als die vorerwähnte Anwendung von
Proportionalmaassstäben in Uebung kam.

Andererseits musste das öftere Betrachten der regelmässigen
Figuren einen geometrisch disponirten Geist von
selbst zum Aufsuchen unbekannter Eigenschaften derselben
reizen, und vielleicht ist der Thaletische Satz von der
Halbirung des Kreises durch einen Durchmesser nichts als
eine aus der Betrachtung jener aegyptischen Zeichnungen
gewonnene Abstraction, und huldigen wir in dieser Beziehung
der Ansicht, dass \Person{Thales} beim Ausspruche des erwähnten,
für uns freilich höchst einfach klingenden Satzes, wahrscheinlich
sagen wollte, nur der Kreis habe die ausgezeichnete
Eigenschaft, von allen durch einen Punkt, den Mittelpunkt,
gehenden Geraden in lauter untereinander gleiche Hälften
getheilt zu werden.

Von besonderer Wichtigkeit scheint uns jedoch der
früher citirte selbstgefällige Ausspruch des \Person{Demokritos} zu
sein, da er uns vor einer ungerechtfertigten Unterschätzung
aegyptischer Constructionsgewandtheit bewahren kann. Bedenklich
in \Person{Demokritos}' Angabe könnte allenfalls jenes
Selbstlob erscheinen, das er sich spendet; wenn es nun
wohl auch schon im Alterthume Männer geben mochte, die
ihre Berühmtheit vorzugsweise und oft nur der Hochschätzung
verdankten, die sie sich selbst und ihren Werken gezollt,
Männer, welche in der Verbreitung des eigenen Lobes so
emsig, so unermüdlich waren, dass sich um sie als die
davon Ueberzeugtesten noch ein Kreis von Gläubigen bildete,
welche den, oft nur auf schwankenden Füssen einhergehenden
Ruhm ihrer Profeten weiter führten, so ist doch die Bedeudung
des Geometers \Person{Demokritos} durch so viele, und verschiedenen
Quellen entspringende Aussprüche beglaubigt,
dass es gewiss Niemandem einfallen wird, seine Autorität
als die eines gründlichen Kenners der Geometrie seiner Zeit
in Zweifel zu ziehen. Wohl sind uns von den geometrischen
Werken des \Person{Demokritos}, und kaum von allen nur die ganz
allgemein klingenden Titel erhalten.

Während uns \Person{Cicero}\footnote{
  \Author{Cicero}, \Title{De finibus bonorum ed malorum}
    I, 6, 20.}
diesen Philosophen als einen
gelehrten, in der Geometrie vollkommen bewanderten Mann
anpreist, theilt uns \Person{Diogenes Laertius}\footnote{
  \Author{Diogenes Laertius}
    IX, 47.}
mit, dass \Person{Demokritos}
»über Geometrie«, »über Zahlen«, »über den Unterschied
des Gnomon oder über die Berührung des Kreises
und der Kugel«, sowie zwei Bücher »über irrationale Linien
und die dichten Dinge« geschrieben habe, Schriften, deren
Titel theilweise uns über ihren Inhalt ganz im Unklaren
lassen. Legen wir den angeführten Zeugnissen Glauben bei,
und es ist kein Grund vorhanden dies nicht zu thuh, so
müssen wir von \Person{Demokritos} als von einem »in der Geometrie
  vollkommenen Manne« voraussetzen, dass er mit den
Errungenschaften des \Person{Pythagoras}, welcher ein Jahrhundert
vor \Person{Demokritos} Aegypten besucht hatte, vollkommen vertraut
war. Gewiss war ihm somit bekannt: die Methode der
»Anlegung der Flächen«, welche wieder die Vertrautheit mit
den Hauptsätzen aus der Theorie der Parallelen und der
Winkel, so wie die Kenntniss der Abhängigkeit der Flächeninhalte
von den ihnen zukommenden Ausmaassen voraussetzt.
Nicht minder bekannt mussten ihm die, dem \Person{Pythagoras}
zugeschriebenen Constructionen der fünf regelmässigen,
sogenannten kosmischen Körper sein, woraus sich weiter
schliessen lässt, dass auch einerseits die Eigenschaften der
Kugel, welcher doch jene Körper eingeschrieben wurden,
und anderseits die Entstehungen der regelmässigen, jene
Körper begrenzenden Vielecke, vor Allem die des Fünfeckes
dem \Person{Demokritos} nicht ungeläufig sein konnten. Die Construction
des Letzteren erheischt wiederum die Kenntniss der
Lehre vom goldenen Schnitt, und diese den Satz vom Quadrate
der Hypothenuse\footnote{
  \Author{Cantor},
    \Title{Vorlesungen über Geschichte der Mathematik},
    I, p. 144--159
    (Leipzig 1880).}.
Hat nun \Person{Demokritos} auch selbst
nichts Neues hinzugefügt, so musste er doch Jenes kennen;
wenn er nun anderseits sagt: »im Construiren hätte ihn
  Niemand, selbst nicht die Harpedonapten der Aegypter übertroffen«,
so dürfen wir hieraus mit Sicherheit schliessen,
dass die geometrischen Kenntnisse der aegyptischen Priester
bedeutend genug gewesen sein mussten, weil sich \Person{Demokritos}
sonst kaum gerade über diese Geometer gesetzt hätte.

Doch verlassen wir für jetzt die Nachrichten des griechischen
Alterthums, welche in der Beurtheilung aegyptischer
Geometrie nur Conjecturen zulassen, und blicken wir nach
directen Denkmalen aegyptischen Ursprungs, aus denen vielleicht
Schlüsse gezogen werden könnten auf Wesen und
Umfang aegyptischer Geometrie.

Das Britische Museum bewahrt eine Papyrusrolle,
welche aus dem Nachlasse des Engländers \Person{A. Henry Rhind}
stammt, die derselbe nebst anderen werthvollen Rollen in
Aegypten käufllich an sich gebracht haben dürfte. Der
erwähnte Papyrus, ein altes Denkmal ägyptischer Mathematik,
ist, wie es scheint, nicht mit vollster Berechtigung als ein
»mathematisches Handbuch« der alten Aegypter bezeichnet
worden\footnote{
  \Author{Eisenlohr},
    \Title{Ein math. Handbuch der alten Aegypter}.
    Leipzig 1877.}.
  Der fragliche Papyrus nennt sich selbst eine
Nachahmung älterer mathematischer Schriften, denn es heisst
in der Einleitung: »Verfasst wurde diese Schrift im Jahre
dreiunddreissig im vierten Monat der Wasserzeit unter König
\Person{Ra-\={a}-us}, Leben gebend nach dem Muster alter Schriften in
den Zeiten des Königs $\ldots$ât vom Schreiber \Person{Aahmes} verfasst
die Schrift.«

Nachdem zuerst Dr. \Person{Birch}\footnote{
  \Author{Birch}, in Lepsius'
    \Title{Zeitschrift für ägypt. Sprache und Alterthum},
    1868, p. 108.}
auf diesen mathematischen
Papyrus durch einen kurzen vorläufigen Bericht aufmerksam
gemacht hatte, wurde der Gegenstand von dem
ausgezeichneten Heidelberger Aegyptologen Dr. \Person{Eisenlohr}
einer eingehenden, höchst schwierigen und zeitraubenden
Untersuchung unterzogen, deren Resultate, was die Uebersetzung
betrifft, unseren gegenwärtigen Betrachtungen zu
Grunde liegen. Bezüglich des Alters des Papyrus hat man
jenes der vorhandenen Abschrift von dem Alter des unbekannten
Originals zu unterscheiden. Nach der von \Person{Eisenlohr}
gegebenen Vervollständigung der in der erwähnten Einleitung
auf das Wort König folgenden Lücke, würde der Herrscher,
unter dessen Regierung das Original entstanden ist, der
König \Person{Ra-en-mat} sein, dessen Regierungszeit \Person{Lepsius}\footnote{
  \Author{Lepsius}, \Title{ägypt. Zeitschrift},
    1871, p. 63.}
auf 2221--2179 v.~Chr.~G. legt. Da ferner der Name
\Person{Ra-\={a}-us} in den bis dahin vorhandenen Königslisten nicht
vorkommt, sah man sich, um die Zeit der Entstehung der
Abschrift wenigstens annähernd angeben zu können, darauf
angewiesen, aus der bekannten Sitte der Aegypter die Eigennamen
der eben herrschenden oder der unmittelbar vorhergegangenen
Regenten zu gebrauchen, Schlüsse zu ziehen.
Und da liess der Name \Person{Aahmes} des Schreibers, sowie auch
die (althieratische) Schrift des Papyrus vermuthen, dass derselbe
um 1700 v.~Chr.~G. entstanden sein dürfte. Die Vermuthung
in Bezug auf das Zeitalter der Abschrift hat sich
nun neueren Forschungen zu Folge vollkommen bestätigt.
Denn \Person{Ra-\={a}-us} wurde als der Hyksoskönig \Person{Apophis}
erkannt, und \Person{Aahmes} dürfte seinen Namen von dem, kurze
Zeit dem Apophis vorhergegangenen Könige \Person{Amasis} entlehnt
haben.

Es erscheint so vollkommen sichergestellt, dass unser
Papyrus aus dem achtzehnten Jahrhundert v.~Chr.~G. stammt.
Die Eingangsworte des Papyrus, welche lauten: »Vorschrift
zu gelangen zur Kenntniss aller dunklen Dinge, aller Geheimnisse,
welche enthalten sind in den Gegenständen«, sowie
die Anordnung des Stoffes in Arithmetik, Planimetrie und
Stereometrie, an welche sich ein, verschiedene Beispiele enthaltender
Theil anschliesst, konnten im ersten Augenblicke
den Gedanken aufkommen lassen, dass wir es vielleicht mit
einem Lehrbuche der Mathematik zu thun haben. Der Umstand
jedoch, dass der Papyrus nur die Zusammenstellung,
allerdings eine in gewissem Grade systematische Zusammenstellung
von Aufgaben nebst ihren Lösungen und den zugehörigen
Proben ist, ohne dass Definitionen oder Lehrsätze
und Beweise vorkommen würden, liess den Papyrus wiederum
als eine Aufgabensammlung, als ein Anleitungsbuch für
Praktiker erscheinen. Man ist noch weiter gegangen, und
stellte die Ansicht auf, der Autor habe bei Abfassung dieser
Schrift vorzüglich an Landleute, welchen die Theorie unzugänglich
war, gedacht. Daraufhin weise nicht nur die Formulirung
des grössten Theiles der Aufgaben, welche Verhältnisse
und Bedürfnisse der Landwirthschaft berücksichtigen, sondern
auch der Schlusssatz des Papyrus, welcher sagt: »Fange das
Ungeziefer und die Mäuse, (vertilge) das verschiedenartige
Unkraut, bitte Gott \Person{Ra} um Wärme, Wind und hohes Wasser«.


Dass wir es nicht mit einem Handbuche, welches dem
damaligen Standpunkte der mathematischen Wissenschaften
in Aegypten entsprechen müsste, zu thun haben, ergibt sich
nicht nur aus dem schon hervorgehobenen Mangel an Definitionen,
Lehrsätzen und Beweisen, ja es fehlt selbst jede Erklärung,
sondern auch aus dem Umstände, dass neben der
richtigen Lösung einzelner Aufgaben die unrichtigen oder
unvollendeten Lösungen derselben oder ähnlicher Aufgaben,
sowie manche Wiederholungen vorkommen. Nur nebenbei
verweisen wir darauf, dass in einem Handbuche unzweifelhaft
wenigstens Anklänge an die erste der Wissenschaften des
Alterthums, an die Astronomie, zu finden sein müssten. Doch
ist von diesem Theile der Mathematik im Papyrus nicht die
geringste Spur zu finden. Aufklärungen über den wahren
Charakter des Originals unseres Papyrus, und eine viele Wahrscheinlichkeit
besitzende Vermuthung über die Entstehung
der uns beschäftigenden Abschrift, verdanken wir dem Scharfsinne
des französischen Aegyptologen Eugène \Person{Revillout}.\footnote{
  \Author{Revillout, Eugène},
    \Title{Revue Egyptologique},
    1881,
    Nr. II et III, p. 304.}

Bei richtiger Erwägung des Umstandes, dass oft auf ein
fehlerlos gelöstes Beispiel, falsche Lösungen ähnlicher Beispiele
folgen, welchen sich dann gewöhnlich eine Reihe von
Uebungsrechnungen anschliesst, Rechnungen die einem Schulpensum
in hohem Grade ähnlich sehen, bei Betrachtung der
Thatsache ferner, wie ein und dasselbe Zahlenbeispiel oft
einigemal und zwar so behandelt wird, dass der Reihe nach
die vorkommenden Zahlenwerthe als die berechneten Resultate
erscheinen, drängt sich uns mit \Person{Eugène Revillout} die
Ueberzeugung auf, dass wir es mit dem Uebungs- oder Aufgabenhefte
eines Zöglings jener Unterrichtshäuser (a·sbo) zu
thun haben, wie deren in so manchem Papyrus Erwähnung
geschieht, und in denen die Schüler, welche später Landwirthe,
Verwalter, Feldmesser oder Constructeure werden
wollten, mit den für ihre künftige Laufbahn notwendigen
Rechnungsoperationen vertraut gemacht wurden. Da dieses
Schulheft selbstverständlich nicht für die Oeffentlichkeit
bestimmt sein konnte, so trägt es auch thatsächlich keinen
Autornamen und keine Jahresangabe; denn, was die in der
Einleitung bezüglich der Zeitperiode, in welcher das Original
entstanden sein sollte, gemachte Erwähnung betrifft, so ist
mehr als wahrscheinlich, dass dieselbe von dem Abschreiber
\Person{Aahmes} herrührt, welcher das Original einige Jahrhunderte
nach seiner Entstehung auffand, und dasselbe, der Mathematik
gewiss ganz unkundig, sammt allen Fehlern abschrieb, zu
diesen noch neue hinzufügend. Nachdem \Person{Aahmes} aus der
Aehnlichkeit der Schriftart des mathematischen Heftes mit
der Schrift anderer ihm bekannten Papyri auf das Alter des
ersteren einen im Ganzen und Grossen nicht unrichtigen
Schluss gezogen haben mochte, so können wir das Ende,
vielleicht auch die Mitte des dritten Jahrtausends v.~Chr.~G.
als jene Zeit betrachten, in welcher das Original der Abschrift
entstanden sein dürfte. Ob \Person{Aahmes} die Abschrift
mit der viel versprechenden Einleitung und der zugleich
praktischen und gottesfürchtigen Schlussregel in der Absicht
versehen hatte, um sie an irgend einen einfachen aegyptischen
Landmann um gutes Geld anzubringen, lassen wir dahingestellt,
und wiederholen nur unsere Uebereinstimmung mit
der Ansicht, dass das Original des Papyrus neben den von
einem Lehrer der Mathematik herrührenden Musterbeispielen,
die sehr oft verunglückten Uebungen eines Schülers enthält,
eines Schülers überdies, der nicht zu den hervorragenden
seiner Glasse gehört haben mochte. Und wie kostbar ist
dennoch dieses altägyptische Schulheft! Wenn wir in aller
Eile eine Skizze seines Inhaltes vorführen sollen, so
müssen wir zunächst die sich auf acht Columnen der oben
erwähnten Einleitung anschliessende Theilung der Zahl 2
durch die Zahlen von 3 bis 99 erwähnen; jeder auftretende
Bruch erscheint in zwei bis vier sogenannte
Stammbrüche, Brüche mit dem Zähler Eins, zerlegt, und sind
die Nenner der letzteren meist gerade Zahlen mit einer
grösseren Divisorenanzahl. Im Anschluss an diese Tabelle
finden wir sechs Beispiele, in denen in Form von Brodvertheilungen
die Division der Zahlen l, 3, 6, 7, 8 und 9 durch
die Zahl 10 gelehrt wird, und es folgt hierauf in 17 Beispielen
die sogenannte Sequem- oder Ergänzungsrechnung,
in welcher es sich darum handelt, Zahlenwerthe zu finden,
die mit gegebenen Werthen durch Addition oder Multiplication
verbunden, andere gegebene Zahlenwerthe liefern. Die nächsten
15 Beispiele gehören der sogenannten \Emphasis{Haurechnung}
an, und finden wir in diesem Abschnitte die Lösungen linearer
Gleichungen mit einer Unbekannten. Zwei weitere, der sogenannten
\Emphasis{Tunnu-} oder Unterschiedsrechnung angehörige
Beispiele belehren uns darüber, dass den alten Aegyptern der
Begriff arithmetischer Reihen nicht fremd war. Es folgen nun
sieben Beispiele über Volumetrie, ebensoviele über Geometrie
und fünf Beispiele über Berechnungen von Pyramiden, also
19 Aufgaben über die wir später noch einige Worte sagen
müssen.

Hieran schliessen sich endlich dreiundzwanzig verschiedenen
Materien entlehnte, Fragen des bürgerlichen Lebens
betreffende Beispiele, wie die Berechnung des Werthes von
Schmuckgegenständen, abermals Vertheilungen von Broden
oder von Getreide, Bestimmung des auf einen Tag entfallenden
Theiles eines Jahresertrages, Berechnungen von Arbeitslöhnen,
Nahrungsmitteln sowie des Futters für Geflügelhöfe.
Einer besonderen Ankündigung werth erscheinen uns in
dieser letzten Abtheilung zwei Beispiele; das eine derselben\footnote{
  \Author{Eisenlohr},
    \Title{Ein math. Handbuch der alten Aegypter}.
    Nr. 64.}
lässt keinen Zweifel darüber aufkommen, dass den alten
Aegyptern die Theorie der arithmetischen Progressionen
vollkommen geläufig war, während wir in dem zweiten\footnote{
  ibid. Nr. 79.}
unter der Aufschrift »eine Leiter« die geometrische Progression
von 7 hoch 1 bis 7 hoch 5 nebst deren Summe vorfinden, wobei
die einzelnen Potenzen eigene Namen: an, Katze, Maus,
Gerste, Maass zu führen scheinen.

Nicht unbemerkt lassen wir endlich die in den Haurechnungen
auftretende Benützung mathematischer Zeichen;
so nach links oder rechts ausschreitender Beine für Addition
und Subtraction, drei horizontale Pfeile für Differenz, sowie
endlich ein besonderes, dem unseren nicht unähnliches
Gleichheitszeichen.

Aus dem geometrischen Theile heben wir zunächst, der
Anordnung des Papyrus nicht folgend, die Flächenberechnungen
von Feldern hervor. Die vorkommenden Beispiele
beziehen sich auf quadratische, rechteckige, kreisrunde und
trapezförmige Felder, deren Flächeninhalte aus ihren Längenmaassen
bestimmt werden. Nachdem in den Aufgaben über
die Berechnung des Fassungsvermögens von Fruchtspeichern
mit quadratischer Grundfläche diese letztere gefunden wird
durch Multiplication der Maasszahl der Seite mit sich selbst,
kann es gar keinem Zweifel unterliegen, dass auch die Fläche
des Rechteckes durch Multiplication der Maasszahlen zweier
zusammenstossender Seiten erhalten wurde, da die Erkenntniss
der Richtigkeit der einen Bestimmungsart, jene der Richtigkeit
der anderen involvirt.

Schon die Betrachtung solcher Proportionalmaassstäbe,
wie wir sie im Grabe \Person{Belzoni} bemerken konnten, hätte
die alten Aegypter, die mit Gleichungen und arithmetischen
Reihen umzugehen wussten, auf die Bestimmung der Fläche
eines Rechteckes aus seinen beiden Seitenlängen mit Nothwendigkeit
führen müssen, und werden wir uns durch den
Umstand, dass im Papyrus der diesbezüglichen Aufgabe eine
zu ihr nicht gehörige Lösung beigefügt ist, durchaus nicht
beirren lassen.

Von hohem Interesse ist die, an mehreren Stellen des
Papyrus vorkommende Methode der Flächenberechnung eines
Kreises, welche zeigt, dass die alten Aegypter mit ziemlicher
Annäherung den Kreis zu quadriren wussten, in der That zu
quadriren, weil sie aus dem Durchmesser eine Länge ableiten,
welche als Seite ein Quadrat liefert, dessen Fläche jener des
Kreises gleichgesetzt wurde. Da sie acht Neuntel des Durchmessers
zur Seite jenes Quadrates machten, so entspricht dies
einem Werthe der Ludolphischen Zahl, welcher dem richtigen
Werthe gegenüber um nicht ganz zwei Hundertstel (um
0,018901) zu hoch gegriffen erscheint; für das dritte Jahrtausend
v.~Chr.~G. und im Vergleiche zu dem Werth $\pi = 3$
der Babylonier, und noch mehr im Vergleiche zu dem Werthe
$\pi = 4$ späterer römischer Geometer, jedenfalls eine nicht zu
unterschätzende Annäherung an den richtigen Werth.

Eine Aufgabe behandelt die Flächenbestimmung des
Dreieckes, wobei das Resultat als das Product zweier Seitenlängen
gefunden wird. Die hier beigefügte Figur\footnote{
  ibid. p. 125.}, welche
in Wirklichkeit ein ungleichseitiges langgestrecktes Dreieck
darstellt, kann ebensowohl als die verfehlte Zeichnung eines
rechtwinkligen wie auch eines gleichschenkligen Dreieckes
betrachtet werden.

Letztere Annahme ist von \Person{Eisenlohr} gemacht und
von \Person{Cantor}\footnote{
  \Author{Cantor},
    \Title{Vorlesungen aus der Geschichte der Mathematik},
    I, p. 49.}
acceptirt worden. Darnach würde sich die
Methode der Dreiecksberechnung der alten Aegypter nur als
eine Näherungsmethode darstellen, und ist auch von beiden
genannten Gelehrten der begangene, in diesem Falle in der
That nicht bedeutende Fehler ermittelt worden.


Wir sind dagegen mit Revillout anderer Meinung.

Mit Rücksicht auf den von uns klar erkannten Charakter
des Originales des Papyrus als eines sehr ungenauen Collegienheftes,
dessen Rechnungen ebensosehr wie die vorkommenden
Zeichnungen von der Mittelmässigkeit seines
Zusammenstellers beredtes Zeugniss ablegen, zweifeln wir
keinen Augenblick, dass die fragliche Figur ein rechtwinkliges
Dreieck vorzustellen hatte. Die mangelhafte Schülerzeichnung
ist durch den Copisten \Person{Aahmes} nur noch
schlechter geworden. Dass ein rechtwinkliges Dreieck gemeint
sein soll, erkennt man übrigens auch aus dem Umstande,
dass in der Figur die Maasszahlen der multiplicirten Seiten
bei den Schenkeln des, vom rechten Winkel nur wenig differirenden
Winkels angesetzt sind, wo doch, wenn es sich hätte
um ein gleichschenkliges Dreieck handeln sollen die Maasszahl
der Schenkel in der Figur gewiss bei beiden Schenkeln zu
finden wäre. Dieselben Gründe bestimmen uns zu der
Annahme, dass die im Papyrus befindliche Flächenberechnung
eines Trapezes eine vollkommen richtige ist, indem es sich
auch hier nur um ein Trapez handeln kann, dessen zwei
parallelen Seiten auf einer der nicht parallelen Seiten senkrecht
stehen. Und warum sollten denn die alten Aegypter
nicht die richtige Art der Flächenberechnung auch beliebiger
Dreiecke gekannt haben?

Konnte man einmal die Fläche eines Rechteckes genau
bestimmen, so musste sich durch einfache Anschauung eines,
durch eine Diagonale zerlegten Rechteckes, von selbst die
Regel zur Flächenbestimmung des rechtwinkligen Dreieckes
ergeben; und wurde nun ein beliebiges schiefwinkliges
Dreieck durch ein Höhenperpendikel in zwei rechtwinklige
zerlegt, so war nichts leichter als die allgemeine Regel zur
Bestimmung der Dreieckfläche aus Basis und Höhe (\Foreign{tepro}
und \Foreign{merit}) zu entwickeln. Dass die Gewinnung des Höhenperpendikels
sowohl bei Constructionen als auch auf dem
Felde den alten Aegyptern nicht unmöglich war, folgt zunächst
aus der grossen Bedeutung der Winkelmaasses (\Foreign{hapt}) für alle
Operationen der praktischen Geometer Aegyptens. Nicht nur,
dass wir in vielen aegyptischen Documenten das Winkelmaass
erwähnt finden, sieht man auch Könige abgebildet, das Winkelmaass
in der Hand, welches von ihnen vielleicht in derselben
Weise durch symbolische Benützung geehrt wurde, wie der
Kaiser von China alljährlich einmal den Pflug zu führen
pflegt. Ein solches Winkelmaass sieht man übrigens auch auf
einem Wandgemälde abgebildet, das eine Schreinerwerkstätte
darstellt,\footnote{
  \Author{Wilkinson},
    \Title{Manners and customs u. s. w.}
    III., p. 144.}
und es unterliegt keinem Zweifel, dass dasselbe
ebensowohl zur Anlegung rechter Winkel als zum Fällen von
Senkrechten benützt worden ist. Aber auch auf freiem Felde
musste den Aegyptern die Construction rechter Winkel geläufig
sein; sowohl die Pyramiden als auch die aegyptischen
Tempel sind vollkommen orientirt, und wurde, wie uns alte
Inschriften\footnote{
  \Author{Brugsch},
    \Title{Ueber Bau und Maasse des Tempels von \Place{Edfu}}
    (\Title{Zeitschrift für ägypt. Sprache u. Alterth.}
    Bd. VIII.)}
belehren, die Orientirung in festlicher Weise
vom Könige unter Beihilfe der Bibliotheksgöttin \Person{Safech} vollzogen,
mit den Worten: »Ich habe gefasst den Holzpflock und
den Stiel des Schlägels, ich halte den Strick gemeinschaftlich
mit der Göttin \Person{Safech}. Mein Blick folgt dem Gange der
Gestirne. Wenn mein Auge an dem Sternbilde des grossen
Bären angekommen ist, und erfüllt ist der mir bestimmte
Zeitabschnitt der Zahl der Uhr, so stelle ich auf die Eckpunkte
Deines Gotteshauses.«

In welchem Maasse bei diesen Operationen die von
\Person{Demokritos} so hochgestellten \Person{Harpedonapten} oder
Seilspanner betheiligt waren, hat \Person{Cantor}\footnote{
  \Author{Cantor},
    \Title{Vorlesungen u. s. w.}
    I, p. 55.} in höchst
scharfsinniger Weise zu beleuchten versucht, und es erscheint
auch uns wahrscheinlich, dass sich die alten Aegypter beim
Construiren rechter Winkel sowie beim Fällen von Senkrechten
auf dem Felde, der Thatsache bedienten, dass der
eine Winkel in einem, die Seitenlängen drei, vier und fünf
besitzenden Dreiecke, ein rechter Winkel sein müsse. Musste
ja doch dieser Satz seit unvordenklichen Zeiten auch den
Chinesen bekannt sein, da wir ihn in der bei ihnen so
berühmten Schrift \Title{Tschiu-p\={\i}} finden, welche mehrere Jahrhunderte
v.~Chr.~G. entstanden, auf den Kaiser \Person{Tsch\={i}u-Kung}
also in das Jahr 1100 v.~Chr.~G. etwa zurückgeführt
wird.\footnote{
  \Author{Éd. Biot}, \Title{Journal Asiatique},
    Paris 1841,
    I. Sem. p. 593.}
Uebrigens konnten directe Messungsversuche an
diagonalen Linien in den Proportionalmaassstäben sowohl zu
dem erwähnten als auch noch zu anderen rechtwinkligen
Dreiecken mit rationalen Seitenlängen geführt haben, und
scheint uns die Möglichkeit nicht ausgeschlossen, dass der
berühmte und berüchtigte Satz des \Person{Pythagoras} über die
Quadrate der Katheten und der Hypothenuse einer eingehenden
Untersuchung solcher Proportionalmaassstäbe entsprungen
ist.

Wenn wir nun einerseits behaupten, dass die alten
Aegypter nicht nur die Fläche des Kreises, des Quadrates,
des Rechteckes, des rechtwinkligen sowie des schiefen Dreieckes,
und unter Zuhilfenahme der Zerlegungen auch die
Flächen beliebiger Polygone theoretisch genau zu bestimmen
im Stande waren, mit Ausnahme der auch für uns eine solche
bildenden Kreisfläche, so muss doch anderseits zugestanden
werden, dass man sich bei praktischen Anwendungen mit
Näherungen begnügte, welche im Laufe der Zeiten so ausarteten,
dass der Gebrauch falscher Regeln ein allgemeiner
wurde.

Am linken Nilufer in der Mitte zwischen \Place{Theben} und
\Place{Assuan} liegt \Place{Edfu}, das alte \Place{Appollinopolis Magna}
mit einem stattlichen Tempelbau aus den Zeiten der Ptolomäer.
Der Tempel, hauptsächlich dem Gotte \Person{Horus} geweiht,
ist mit einer freistehenden Umfassungsmauer umgeben,\footnote{
  \Author{Lepsius},
    \Title{Ueber eine hieroglyphische Inschrift am Tempel von \Place{Edfu}}.
    \Title{Abhandlung d. Acad. d. Wiss. in Berlin},
    1855, p. 69.}
deren Ostseite zwischen dem Brunnenthore und dem östlichen
Pylonflügel eine Inschrift trägt, welche uns auf acht
Feldern und in hundertvierundsechzig Columnen\footnote{
  \Author{Brugsch}, \Title{Thesaurus III},\label{br}
    Leipzig 1884.} eine
Schenkungsurkunde des Königs \Person{Ptolomäus XI. Alexander
  I.} (mit dem Beinamen \Person{Philometor}) bekannt gibt. Das
Geschenk, welches hier \Person{Horus} und den übrigen Göttern von
\Place{Edfu} verliehen wird, besteht aus einer Anzahl von meist
viereckigen Aeckern, deren vier Seitenlängen nebst Flächeninhalten
angegeben erscheinen.

Da jeder der vorkommenden Flächeninhalte identisch
ist mit dem Producte der arithmetischen Mittel der beiden
Gegenseitenpaare, so wurde nach \Person{Lepsius} die Vermuthung
aufgestellt, die alten Aegypter hätten, um Vierecke bei der
Flächenbestimmung annähernd wie Rechtecke behandeln zu
können, den Unterschied der Gegenseiten dadurch auszugleichen
gesucht, dass sie die arithmetischen Mittel derselben
in Rechnung zogen.

Bei sehr vielen der in der \Place{Edfu}er Schenkungsurkunde
vorkommenden Vierecke ist der Unterschied je zweier Gegenseiten
entweder Null oder verhältnissmässig so klein, dass
man den betreffenden Vierecken eine vom Rechtecke wenig
verschiedene Gestalt beilegen kann, und die erhaltenen Resultate
somit eine ziemliche Annäherung an den richtigen
Flächenwerth darstellen dürften, nach dem man mit Rücksicht
auf die bei \Person{Sesostris} bemerkte Eintheilung des Landes
in Rechtecke voraussetzen darf, gerade diese oder eine ihr
zunächst kommende Form der Felder sei die auch damals
schon beliebte gewesen.

Doch kommen auch Vierecke vor, wo der Längenunterschied
der Gegenseiten ein bemerkenswerther ist; ja es werden
auch Dreiecke als Vierecke mit einer verschwindenden Seite
behandelt, so dass der begangene Fehler in manchen Fällen
ein nicht unbedeutender ist.

Nur nebenbei bemerken wir, dass man dieselbe unrichtige
Flächenformel für das Viereck erhält, wenn man dasselbe
zunächst durch eine Diagonale in zwei Dreiecke zerlegt,
auf jedes dieser Dreiecke die unrichtige Flächenformel,
die den Inhalt als das halbe Product der beiden Seiten liefert,
anwendet, die beiden so erhaltenen Dreiecksflächen addirt
und dann aus dieser Summe und jener, welche man bei dem
ähnlichen Vorgange durch Zerlegung mittelst der zweiten
Diagonale erhält, das arithmetische Mittel construirt.

Nimmt man mit \Person{Eisenlohr} und \Person{Cantor} an, dass die
Aegypter die Dreiecksfläche wirklich dem halben Producte
zweier Seiten gleichsetzten, so steht man vor der Frage,
warum nicht in derselben Art die Flächen der in der \Place{Edfu}er
Schenkungsurkunde auftretenden Dreiecke bestimmt erscheinen?

Uebrigens wolle man sich darüber nicht wundern, dass
es überhaupt möglich war, die Flächenberechnungen im
praktischen Leben nach einer so falschen Methode durchzuführen.
Wissen wir doch, dass im Alterthume, zur Zeit
\Person{Platon}s, einer der gebildetsten Männer, einer der
hervorragendsten Geschichtschreiber, dass \Person{Thukydides}\footnote{
  \Author{Thukydides},
    ed. \Editor{Rothe},
    VI. 1.} in seiner
Unkenntniss der Beziehung zwischen Flächeninhalt und Umfang,
die Fläche einer Insel nach der zu ihrer Umschiffung
nothwendigen Zeit zu bestimmen suchte; in der Geometrie
\Person{Gerbert}'s,\footnote{
  ed. \Editor{Olleris},
    Cap. LXX. p. 460.}
des nachmaligen Papstes \Person{Silvester II.} 
finden wir, 1000 Jahre nach Chr.~G., die Fläche eines gleichschenkligen
Dreieckes durch Multiplication des Schenkels
mit der halben Basis berechnet, wo doch schon \Person{Hero von
Alexandrien}\footnote{
\Title{\Emphasis{Heronis Alexandrini} geometricorum et stereometricorum reliquiae}
  (ed. \Editor{Hultsch},
  Berlin 1864).}
1100 Jahre früher die richtige Formel
für diese Berechnung kennt.

Wir berühren diese Thatsachen, und könnten noch eine
ganze Reihe ähnlicher Beispiele anführen, nur um zu zeigen,
wie übereilt es wäre, aus den oft nur schwache Annäherungen
liefernden Berechnungen der \Place{Edfu}er Schenkungsurkunde
schliessen zu wollen, die richtigen Methoden seien
den in die Wissenschaften eingeweihten aegyptischen Priestern
nicht bekannt gewesen.

Doch zurück zum Papyrus \Emphasis{Rhind}.

Wir übergehen die Inhaltsbestimmungen von Fruchthäusern,
bei denen der Inhalt durch Multiplication einer
Fläche mit einer Länge bestimmt wird, weil wir es für
müssig halten, Erörterungen darüber anzustellen, welche
Flächen und Längen hiebei gemeint sind, so lange uns über
die Form jener Fruchthäuser oder Speicher nichts bekannt
ist.

Dagegen erwecken die im Papyrus vorkommenden Pyramiden-Berechnungen
das höchste Interesse, besonders nach
den glänzenden Untersuchungen, welchen \Person{Revillout}
diesen Gegenstand unterzogen hat, und deren Resultate wir,
entgegen der von \Person{Eisenlohr} ausgesprochenen und auch
von \Person{Lepsius}\footnote{
  \Author{Lepsius},
    \Title{Ueber die 6palmige grosse Elle von 7 kleinen Palmen Länge
    in dem »math. Handbuche« von Eisenlohr}.
  (\Title{Zeitschrift f. äg. Sp.}
  1884. 1. Heft.)}
acceptirten Ansicht als solche betrachten,
welche in einfacher und natürlicher Weise die sogenannte
\Emphasis{Seket}-Rechnung der alten Aegypter beleuchten.

Es wird in diesen Rechnungen die Böschung der Seitenflächen
einer quadratischen Pyramide dadurch fixirt, dass
jener Theil der Länge eines der beiden gleichlangen Schenkel
des Winkelmaasses berechnet wird, der sich zur Länge des
anderen Schenkels so verhält, wie die halbe Länge der Basisseite
der quadratischen Pyramide zur Höhe derselben.


Zu dem Behufe war der eine der beiden Schenkel des
Winkelmaasses in eine gewisse Anzahl gleich grosser Theile
getheilt, während der andere Schenkel, der Pyramidenhöhe
entsprechend, und als Einheit betrachtet, ungetheilt blieb.

Um nun den sogenannten \Emphasis{Seket} zu bestimmen, wurde
die halbe Länge der Basisseite durch die Pyramidenhöhe
dividirt und mit dem erhaltenen Quotienten die Anzahl der
Theile des horizontalen, getheilten Schenkels des Winkelmaasses
multiplicirt.

Es war somit der Seket (welcher in derselben Art für
einen geraden Kreiskegel aus dem Durchmesser der Basis
und der Höhe bestimmt erscheint) als Verhältniss aufgefasst,
die goniometrische Cotangente des Neigungswinkels der Seitenfläche
der Pyramide, respective der Kegelkante zur Basis.

Wenn wir selbstverständlich weit davon entfernt sind,
hierin vielleicht Anfänge der Trigonometrie sehen zu wollen,
so erkennen wir doch anderseits, dass den alten Aegyptern
auch die Lehre proportionaler Linien, wenigstens in ihren
Anwendungen, bekannt gewesen sein musste, und erscheint
uns auch der am Eingange erwähnte Ausspruch über die
dem Milesier \Person{Thales} zugeschriebene Höhenmessung der
Pyramiden als ein ganz glaubwürdiger, wenn wir sehen, wie
im Papyrus von den drei Werthen: Basis, Höhe, Seket, jeder
aus den beiden anderen berechnet erscheint.

Fassen wir nun die Ergebnisse unserer Betrachtungen
zusammen, so müssen wir aus der quellenmässig erwiesenen
grossen Bewunderung, welche die ausgesprochen geometrisch
hochentwickelten Griechen den aegyptischen Geometern rückhaltlos
zollten, wir müssen aus der unanfechtbaren Thatsache,
dass griechische Geometer den Grund zu ihren Kenntnissen
und Entdeckungen in Aegypten suchten und fanden,
wir müssen im Hinblicke auf das, aus der nun vollends
entzifferten\footref{br} \Place{Edfu}er Schenkungsurkunde sich mit Sicherheit
ergebende ausgebreitete und fest organisirte Katasterwesen
der alten Aegypter, welches zugleich mit den zahlreichen,
dem öffentlichen Leben dienenden Land- und
Wasserbauten auf eine verhältnissmässig bedeutend entwickelte
Vermessungskunde hinweist, wir müssen endlich
aus dem von uns besprochenen Papyrus, der sich als eine
ungenaue Abschrift eines mangelhaften, aus dem dritten Jahrtausend
vor Chr.~G. stammenden, mathematischen Collegien- oder
Aufgabenheftes erweist, und aus dessen Vorhandensein
sich fast mit Gewissheit auf damals existirende, neben den
Regeln auch ihre Ableitungen enthaltende Lehrbücher
schliessen lässt, wir können und müssen aus allen diesen
Umständen den allgemeinen Schluss ziehen, dass bereits drei
Jahrtausende vor unserer Zeitrechnung sowohl die arithmetischen,
als auch die geometrischen Kenntnisse der Aegypter,
einen für dieses Zeitalter bedeutenden Grad der Entwicklung
besassen.

Insbesondere können wir in jenen fernen Zeiten eine
staunenswerth weitgehende Annäherung bei der Berechnung
der Kreisfläche beobachten, wir finden mit vollständiger
Sicherheit richtige Flächenbestimmungen des Quadrates,
Rechteckes und des rechtwinkligen Dreieckes; höchst wahrscheinlich
auch richtige Bestimmungen der Flächen schiefwinkliger
Dreiecke und Vierecke, welche im praktischen
Leben durch leichter zu handhabende Annäherungsformeln
ersetzt wurden; wir sehen Bestimmungen des Rauminhaltes
durch ihre Dimensionen gegebener Körper und erkennen die
Anfänge der Aehnlichkeitslehre.

Was das geometrische Zeichnen betrifft, so kennen wir
schon die Construction der früher beobachteten regelmässigen
Figuren und dürfen weiter vermuthen, dass die Anlegung
rechter Winkel und das Fällen von Senkrechten sowohl
mittelst des Winkelmaasses als auch mittelst rationaler rechtwinkliger
Dreiecke bekannt, und die Zerlegung gegebener
Flächen behufs ihrer Inhaltbestimmung in allgemeiner Verwendung
war.

Gewiss werden auch theoretische Resultate bekannt
gewesen sein; so die Hälftung des Kreises durch seinen
Durchmesser, die sich aus der besprochenen Seketrechnung
von selbst ergebende Winkelgleichheit an der Basis gleichschenkliger
Dreiecke und gleichseitiger quadratischer Pyramiden,
und wohl noch manches Andere.

Möge es gelingen, durch Auffindung neuer, sowie durch
Entzifferung der, noch ihrer Erklärung harrenden Denkmale
und Schriften, von welchen letzteren, Dank der hohen Munificenz
des Erlauchten Curators unserer Akademie, auch Wien
eine imposante Zahl aufweisen kann, möge es so gelingen
noch weitere Anhaltspunkte für die Kenntniss der mathematischen
Thätigkeit des uns bekannten ältesten Culturvolkes,
der Aegypter zu gewinnen!

Diesen unseren Wunsch theilen gewiss Alle, denen die
Erforschung der Culturgeschichte des menschlichen Geschlechtes
nicht ohne Wichtigkeit erscheint!

\clearpage
\begin{center}
\textbf{Unvollständige bibliografische}\\
\textbf{Angaben im Text}\footnote{Unvollständige
bibliografische Angaben haben wir, wo möglich, ergänzt und hier zusammengestellt. RS}
\end{center}
\bigskip

\begin{itemize}
\item \Author{Plato} (ed. \Editor{Friedrich Ast}), \Title{Platonis Phaedrus},
Leipzig, 1810/1830.
\item \Author{Heron, Alexandrinus} et \Author{Didymus, Alexandrinus} (ed.
\Editor{Friedrich Hultsch}), \Title{Heronis Alexandrini Geometricorum et Stereometricorum reliquiae : accedunt Didymi Alexandrini Mensurae marmorum et anonymi variae collectiones ex Herone, Euclide, Gemino, Proclo, Anatolio aliisque},
Berlin, Weidmann, 1864.
\item \Author{Strabo} (ed. \Editor{August Meineke}), \Title{Strabonis Geographia}, Leipzig, Teubner, 1852--77.
\item \Author{Eudemus} von Rhodos (ed. \Editor{Leonhard von Spengel}), \Title{Eudemi Rhodii Peripatetici fragmenta quae supersunt}, Berlin (S.~Calvary), London (Williams \& Norgate), 1866.
\item \Author{Proklos} (ed. \Editor{Gottfried Friedlein}), \Title{Procli Diadochi in primum Euclidis elementorum librum commentarii}. Leipzig 1873.
\item \Author{Montucla, Jean Etienne} et \Author{Lalande, J}, \Title{Histoire des mathématiques, dans laquelle on rend compte de leurs progrès depuis leur origine jusqu'à nos jours $\ldots$}, Paris, 1799--1802.
\item \Author{Bretschneider, Carl Anton}, \Title{Die Geometrie und die Geometer vor Euklides : ein historischer Versuch}, Leipzig, Teubner, 1870.
\item \Author{Cantor, Moritz}, \Title{Vorlesungen über Geschichte der Mathematik}, Leipzig, Teubner, 1880.
\item \Author{Clemens, Alexandrinus} (ed. \Editor{John Potter}, \Title{Ta Heuriskomena}, Oxford, 1715.
\item \Author{Theon} von Smyrna (ed. \Editor{Thomas H. Martin}), \Title{	
Liber de astronomia cum Sereni fragmente}, Paris, 1849.
\item \Author{Avennes, Prisse d'}, \Title{Histoire de l'art égyptien d'après les monuments}, Paris, 1879.
\item \Author{Brugsch, Heinrich}, \Title{Thesaurus Inscriptionum Aegyptiacarum}, Leipzig, Hinrichs, 1884.
\item \Author{Sylvester, Papa II} (Gerbert d'Aurillac) (ed. \Editor{Alexandre Olleris)}, \Title{\OE{}uvres de Gerbert, Pape sous le nom de Sylvestre II}, Clermont-Ferrand, Thibauld, 1867.
\end{itemize}
\end{document}